\documentclass[11pt,a4paper]{amsart}
\usepackage[T1]{fontenc}
\usepackage{lmodern}
\usepackage{amsmath,amssymb,mathtools}
\usepackage[margin=28mm]{geometry}
\usepackage{microtype}
\usepackage{xcolor}
\usepackage[colorlinks=true,linkcolor=blue!45!black,citecolor=blue!45!black,urlcolor=blue!45!black]{hyperref}
\usepackage{enumitem}
\usepackage{needspace}
\hypersetup{pdftitle={On local embeddability into C*(F-infinity)},pdfauthor={Jean Roydor},pdfsubject={Local embeddability into the full free group C*-algebra}}
\newtheorem{theorem}{Theorem}[section]
\newtheorem{proposition}[theorem]{Proposition}
\newtheorem{lemma}[theorem]{Lemma}
\newtheorem{corollary}[theorem]{Corollary}
\theoremstyle{definition}
\newtheorem{definition}[theorem]{Definition}
\theoremstyle{remark}
\newtheorem{remark}[theorem]{Remark}
\numberwithin{equation}{section}
\newcommand{\Z}{\mathbb Z}
\newcommand{\N}{\mathbb N}
\newcommand{\F}{\mathbb F}

\newcommand{\U}{\mathcal U}
\newcommand{\id}{\mathrm{id}}
\newcommand{\SL}{\mathrm{SL}}
\newcommand{\cb}{\mathrm{cb}}
\newcommand{\tmin}{\otimes_{\min}}
\newcommand{\tmax}{\otimes_{\max}}
\newcommand{\tM}{\otimes_M}
\DeclareMathOperator{\tr}{tr}
\DeclareMathOperator{\spanop}{span}
\DeclareMathOperator{\Ree}{Re}
\setlist[enumerate]{label=\textup{(\roman*)},leftmargin=*,itemsep=3pt}
\title{On local embeddability into $C^*(\F_\infty)$}
\author{Jean Roydor}
\address{Institut de Math\'ematiques de Jussieu, Sorbonne Universit\'e}
\email{jroydor@imj-prg.fr}
\date{}
\begin{document}
\begin{abstract}
We study local embeddability of $C^*$-algebras into $C^*(\F_\infty)$
as an analogue of QWEP in the LLP/WEP picture. We characterize this property by an analog of CP-stability,
an approximation property for finite-dimensional operator systems.
We show that local embeddability,
together with an amenable trace, implies a tracial local lifting property.
Combining this observation with the cohomological obstructions of Ioana,
Spaas and Wiersma, we prove that the full group $C^*$-algebra
$C^*(\SL_3(\Z))$ does not locally embed into
$C^*(\F_\infty)$. We obtain the same conclusion for the full
$C^*$-algebras of two hyperlinear groups constructed by Thom.
To our knowledge, the only previously known obstruction to local embeddability
into $C^*(\F_\infty)$ was due to Junge and Pisier.
\end{abstract}
\maketitle

\section{Introduction}

Kirchberg's work \cite{Kirchberg} revealed a close duality between the
local lifting property (LLP) and the weak expectation property (WEP),
expressed through tensor products of $C^*$-algebras. The quotient weak
expectation property (QWEP) extends the WEP side of this picture. In this
paper, we study local embeddability into $C^*(\F_\infty)$ as its counterpart
on the LLP side.

%Let $\F_\infty$ denote the free group on countably many generators.
We first recall the operator-space invariant that quantifies local
embeddability.

For an operator space $E$, a linear map $T:E\to F$ is completely bounded if
$\|T\|_{\cb}=\sup_{r\ge1}\|\id_{M_r}\otimes T\|<\infty$.
The completely bounded Banach--Mazur distance between finite-dimensional
operator spaces is
\[
d_{\cb}(E,F)=\inf_T\|T\|_{\cb}\|T^{-1}\|_{\cb},
\]
where $T$ runs over linear cb-isomorphisms from $E$ onto $F$.

For a finite-dimensional operator space $E$, set
\[
d_f(E)=\inf\{d_{\cb}(E,F):F\subset C^*(\F_\infty),\ \dim F=\dim E\}.
\]
For an operator space $X$, set
\[
d_f(X)=\sup\{d_f(E):E\subset X,\ \dim E<\infty\}.
\]
In particular, this defines $d_f(A)$ for a $C^*$-algebra $A$, with its
canonical operator-space structure.\\
We say that $A$ \emph{locally embeds into $C^*(\F_\infty)$} if $d_f(A)=1$.
Equivalently, for every finite-dimensional $E\subset A$ and
every $\varepsilon>0$, there is an injective linear map $T:E\to C^*(\F_\infty)$
such that
\begin{equation}\label{eq:localembedding}
\|T\|_{\cb}\|T^{-1}:T(E)\to E\|_{\cb}<1+\varepsilon.
\end{equation} 
%Neither multiplicativity nor positivity is required of $T$.
Exact $C^*$-algebras locally embed into $C^*(\F_\infty)$. Moreover,
every separable nuclear $C^*$-algebra admits a completely isometric
linear embedding into $C^*(\F_\infty)$. For details, see Section~20 in
\cite{Pisier}; see also \cite{Harcharras}.

%: in the unital case, the
%Choi--Effros lifting theorem \cite{CE} gives a u.c.p. section of a
%quotient map from $C^*(\F_\infty)$ onto the algebra; the general case
%follows by unitization.
Every $C^*$-algebra with the local lifting property (LLP) locally embeds
into $C^*(\F_\infty)$. The converse fails: the reduced group algebra $C_r^*(\F_2)$ is exact
and therefore locally embeds into $C^*(\F_\infty)$, but it does not have
the LLP \cite{Kirchberg}.

To our knowledge, the only previously known obstruction to local
embeddability into $C^*(\F_\infty)$ is due to Junge and Pisier, who proved
in 1995 that $\mathbb B(\ell_2)$ does not locally embed into $C^*(\F_\infty)$
\cite{JP}. Their work also provides a tensor-norm characterization
of local embeddability, recalled in Proposition~\ref{prop:df-dichotomy}.
It implies that the invariant $d_f(A)$ can take only
the values $1$ and $\infty$ for a $C^*$-algebra $A$.

In Section~\ref{sec:characterisations}, we recall characterizations of local embeddability into $C^*(\F_\infty)$
by tensor products with WEP algebras and by the local liftability of
WEP maps (Theorem~\ref{thm:df-WEP-maps}). Comparison with Pisier's QWEP
characterization (Theorem~\ref{thm:QWEP-Pisier}) makes the parallel more explicit.
We also establish permanence properties (Proposition~\ref{prop:df-permanence}).
Section~\ref{sec:ACP} gives an equivalent formulation for unital
$C^*$-algebras in terms of ambient CP-stability (ACP-stability).
The approximation condition allows the controlling finite-dimensional
operator system to lie in an ambient $\mathbb B(H)$.

In Section~\ref{sec:tracial-lifting}, we obtain new obstructions for
full group $C^*$-algebras. We prove that $C^*(\SL_3(\Z))$ and
$C^*(\Z^2\rtimes\SL_2(\Z))$ do not locally embed into $C^*(\F_\infty)$
(Theorem~\ref{thm:SL3-tracial-lifting}). The key observation is that
local embeddability, together with an amenable trace, forces local
lifting for trace-preserving representations into finite von Neumann
algebras (Proposition~\ref{prop:trace-LLP}). This is a reminiscence of the well-known fact that local embeddability, together with WEP imply the LLP. We can incorporate this key observation to the
projective-representation obstruction of Ioana, Spaas and Wiersma
\cite{ISW}. We also treat Thom's two examples of hyperlinear groups
without the factorization property \cite{Thom}, in
Corollary~\ref{cor:Thom-second} and Section~\ref{sec:Thom-first}.
For the first example, the argument uses an amenable trace other than
the canonical group trace.

%All groups are discrete, all representations of unital $C^*$-algebras are
%unital, and Hilbert-space inner products are linear in the second variable.
%The notation $C^*(G)$ always means the full group $C^*$-algebra; its canonical
%unitaries are denoted by $u_g$, $g\in G$.

\Needspace{6\baselineskip}
We recall the tensor norm that detects local embeddability into $C^*(\F_\infty)$, see \cite{Pisier} for details. If
$E\subset \mathbb B(H)$ and $F\subset \mathbb B(K)$ are completely isometric embeddings, the
Junge--Pisier norm $M$ on $E\odot F$ is the norm induced by
$\mathbb B(H)\tmax \mathbb B(K)$. It is independent of the embeddings, is injective for
complete isometries, and satisfies
\begin{equation}\label{eq:Mfunctorial}
\|T\otimes S:E\tM F\longrightarrow E'\tM F'\|_{\cb}
\le\|T\|_{\cb}\|S\|_{\cb}.
\end{equation}
For $C^*$-algebras, $M$ is a $C^*$-tensor norm and consequently 
$\|\cdot\|_{\min}\le\|\cdot\|_M\le\|\cdot\|_{\max}$.\\
The next proposition is probably well-known to specialists, but we wanted to recall this dichotomy
\Needspace{6\baselineskip}
\begin{proposition}\label{prop:df-dichotomy}
For a $C^*$-algebra $A$, the following assertions are equivalent.
\begin{enumerate}
\item $d_f(A)<\infty$.
\item $d_f(A)=1$, $i.e.$, $A$ locally embeds into $C^*(\F_\infty)$.
\item $A\tM \mathbb B(\ell_2)=A\tmin \mathbb B(\ell_2)$ isometrically via the identity
on $A\odot \mathbb B(\ell_2)$.
\end{enumerate}
\end{proposition}

\begin{proof}
Assume \textup{(i)} and put $c=d_f(A)$. By \cite[Theorem~4.5]{JP},
$\|x\|_M\le c\|x\|_{\min}$ for every $x\in A\odot \mathbb B(\ell_2)$.
For $z\in A\odot \mathbb B(\ell_2)$ and $n\ge1$, apply this estimate to
$(z^*z)^n$. Since both norms are $C^*$-norms,
\[
\|z\|_M^{2n}=\|(z^*z)^n\|_M
\le c\,\|(z^*z)^n\|_{\min}
=c\,\|z\|_{\min}^{2n}.
\]
Taking $2n$th roots and letting $n\to\infty$ yields
$\|z\|_M\le\|z\|_{\min}$. The reverse inequality always holds, proving
\textup{(iii)}.

Conversely, \textup{(iii)} implies \textup{(ii)} by the same result with $c=1$. Finally, \textup{(ii)} implies \textup{(i)}.
\end{proof}

\begin{remark}[Tensor characterization of QWEP]\label{rem:qwep-L}
A unital $C^*$-algebra $Q$ is QWEP if and only if
\[
Q\otimes_L C^*(\F_\infty)=Q \tmax C^*(\F_\infty)
\]
isometrically on $Q \odot C^*(\F_\infty)$. For quotients
$A=C^*(\F_\alpha)/I$ and $B=C^*(\F_\beta)/J$, where $\F_\alpha$ and
$\F_\beta$ are free groups on $\alpha$ and $\beta$ generators, respectively,
and $I,J$ are closed ideals, the tensor norm $L$ on $A\odot B$ is induced by
the quotient
\[
\frac{C^*(\F_\alpha)\tmin C^*(\F_\beta)}
{I\tmin C^*(\F_\beta)+C^*(\F_\alpha)\tmin J}.
\]
This norm is independent of the free presentations.
Here $L$ is Pisier's notation for the tensor norm denoted by
$\mathrm{lr}$ in \cite{Kirchberg}.
See also \cite[\S7.6 and Theorem~9.67]{Pisier}.
\end{remark}

\section{Characterisations and permanence properties}
\label{sec:characterisations}

We first give characterisations of local embeddability and QWEP in terms
of WEP and LLP maps.

Let $i_B:B\to B^{**}$ denote the canonical
embedding of a $C^*$-algebra, and denote the decomposable norm by
$\|\cdot\|_{\mathrm{dec}}$. For a linear map $u:X\to B$ from an
operator space, the following terminology is that of
\cite[Definitions~9.55 and~9.58]{Pisier}. We call $u$ a \emph{WEP map}
if $\id_{C^*(\F_\infty)}\otimes u:C^*(\F_\infty)\tmin X\to C^*(\F_\infty)\tmax B$
is contractive, and an \emph{LLP map} if
$\id_{\mathbb B(\ell_2)}\otimes u:\mathbb B(\ell_2)\tmin X\to\mathbb B(\ell_2)\tmax B$ is contractive.

Given a quotient map $q:D\twoheadrightarrow D/I$, a linear map
$u:X\to D/I$ from an operator space is \emph{locally $1$-liftable through $q$}
if, for every finite-dimensional subspace $E\subset X$, there is a
complete contraction $s_E:E\to D$ such that $qs_E=u|_E$.

\begin{theorem}[Characterization of local embeddability into $C^*(\F_\infty)$]
\label{thm:df-WEP-maps}
Let $A$ be a unital $C^*$-algebra. The following are equivalent.
\begin{enumerate}
\item $A$ locally embeds into $C^*(\F_\infty)$.
\item For all unital $C^*$-algebras $W_1,W_2$ with WEP,
every u.c.p. map $u:A\to W_1$ induces a contraction
\[
\id_{W_2}\otimes u:W_2\tmin A\longrightarrow W_2\tmax W_1.
\]
\item[\textup{(ii$'$)}] Idem with $W_1=W_2=\mathbb B(\ell_2)$.
\item There exist a Hilbert space $H$ and a faithful unital
$*$-representation $i:A\hookrightarrow \mathbb B(H)$ such that $i$ is an LLP map.
\item For every quotient map $q:D\twoheadrightarrow D/I$ of
$C^*$-algebras, every WEP map $u:A\to D/I$ is locally $1$-liftable
through $q$.
\end{enumerate}
\end{theorem}

Now, to show that local embeddability into $C^*(\F_\infty)$ plays the analogous role of the QWEP in the LLP/WEP picture, we recall Pisier's QWEP criterion
\cite[Theorem~9.67]{Pisier} in its u.c.p. formulation and also include the extension property in
\textup{(iv)} in order to make the these two theorems symmetric to each other.

\begin{theorem}[Characterization of QWEP]
\label{thm:QWEP-Pisier}
Let $Q$ be a unital $C^*$-algebra. The following conditions are equivalent.
\begin{enumerate}
\item $Q$ is QWEP.
\item For all unital $C^*$-algebras $L_1,L_2$ with LLP,
every u.c.p. map $u:L_1\to Q$ induces a contraction
\[
\id_{L_2}\otimes u:L_2\tmin L_1\longrightarrow L_2\tmax Q.
\]
\item[\textup{(ii$'$)}] Idem with $L_1=L_2=C^*(\F_\infty)$.
\item[\textup{(iii)}] There exist a free group $\mathbb F$ and a unital
surjective $*$-homomorphism $q:C^*(\mathbb F)\to Q$ which is a WEP map.
\item[\textup{(iv)}] For every inclusion of operator spaces $X\subset Y$
and every LLP map $u:X\to Q$, there is a complete contraction
$\widetilde u:Y\to Q^{**}$ such that
\[
\widetilde u|_X=i_Q\circ u.
\]
\end{enumerate}
\end{theorem}

Comparing assertions $(iv)$ in both theorems above, we see that local liftability is a weak form of liftability just as bidual valued extension is weak form of extension property.

\begin{lemma}\label{lem:WEP-map-M}
Let $W$ be a $C^*$-algebra with WEP and let $u:X\to B$ be a WEP map.
Then $\id_W\otimes u:W\tM X\longrightarrow W\tmax B$ is contractive.
\end{lemma}

\begin{proof}
By \cite[Corollary~9.56]{Pisier}, the map $i_Bu:X\to B^{**}$ factors as
$wv$, where $v:X\to \mathbb B(K)$ is completely bounded, $w:\mathbb B(K)\to B^{**}$
is decomposable, and $\|v\|_{\cb}\|w\|_{\mathrm{dec}}\le1$.
In a faithful representation $W\subset \mathbb B(H)$, WEP makes this inclusion
isometric after maximal tensoring with any $C^*$-algebra
\cite[Theorem~9.22]{Pisier}. Hence $W\tM \mathbb B(K)=W\tmax \mathbb B(K)$.
The canonical inclusion $B\to B^{**}$ is also isometric after maximal
tensoring \cite[Corollary~7.27]{Pisier}. Thus, for $z\in W\odot X$,
\[
\begin{aligned}
\|(\id_W\otimes u)z\|_{W\tmax B}
 &=\|(\id_W\otimes w)(\id_W\otimes v)z\|_{W\tmax B^{**}}\\
 &\le \|w\|_{\mathrm{dec}}\,\|(\id_W\otimes v)z\|_{W\tM \mathbb B(K)}
 \le \|z\|_{W\tM X},
\end{aligned}
\]
by maximal tensorization of decomposable maps and \eqref{eq:Mfunctorial}.
\end{proof}

\begin{proof}[Proof of Theorem~\ref{thm:df-WEP-maps}]
Assume \textup{(i)}. By \cite[Theorem~4.5]{JP}, $W\tM A=W\tmin A$
for every $C^*$-algebra $W$. In \textup{(ii)}, $u$ is a WEP map:
it is completely contractive and $C^*(\F_\infty)\tmin W_1=C^*(\F_\infty)\tmax W_1$
by WEP of $W_1$. Apply Lemma~\ref{lem:WEP-map-M} with $W=W_2$ to get
\textup{(ii)}. Taking $W_1=W_2=\mathbb B(\ell_2)$ gives \textup{(ii$'$)}.

Assume \textup{(ii$'$)} and fix $z=\sum_{j=1}^m b_j\otimes a_j\in(\mathbb B(\ell_2)\odot A)\setminus\{0\}$.
The unital algebra $A_0=C^*(1,a_1,\ldots,a_m)$ is separable. Choose a
faithful unital representation $\pi:A_0\to\mathbb B(\ell_2)$ and extend it by Arveson's
theorem to a u.c.p. map $u:A\to\mathbb B(\ell_2)$. Since $u|_{A_0}=\pi$ is completely
isometric, injectivity of $M$ and the identity $\mathbb B(\ell_2)\tM\mathbb B(\ell_2)=\mathbb B(\ell_2)\tmax\mathbb B(\ell_2)$ give
\[
\|z\|_{\mathbb B(\ell_2)\tM A}
 =\|(\id_{\mathbb B(\ell_2)}\otimes u)z\|_{\mathbb B(\ell_2)\tmax\mathbb B(\ell_2)}
 \le \|z\|_{\mathbb B(\ell_2)\tmin A}.
\]
The inequality is \textup{(ii$'$)}. Proposition~\ref{prop:df-dichotomy} now
gives \textup{(i)}.

For any faithful unital representation $i:A\hookrightarrow \mathbb B(H)$,
the definition of $M$ gives
\[
\|(\id_{\mathbb B(\ell_2)}\otimes i)z\|_{\mathbb B(\ell_2)\tmax \mathbb B(H)}=\|z\|_{\mathbb B(\ell_2)\tM A}
\qquad(z\in\mathbb B(\ell_2)\odot A).
\]
Thus $i$ is an LLP map exactly when $\|z\|_{\mathbb B(\ell_2)\tM A}\le\|z\|_{\min}$
for all such $z$. By Proposition~\ref{prop:df-dichotomy}, this is equivalent
to \textup{(i)}. This proves \textup{(i)}$\Leftrightarrow$\textup{(iii)}
and shows that every faithful unital representation works if one does.

To see that \textup{(i)} implies \textup{(iv)}, put $B=D/I$.
Lemma~\ref{lem:WEP-map-M} with $W=\mathbb B(\ell_2)$ and \cite[Theorem~4.5]{JP}
give a contraction $\id_{\mathbb B(\ell_2)}\otimes u:\mathbb B(\ell_2)\tmin A\to\mathbb B(\ell_2)\tmax B$.
Maximal tensoring preserves quotients, so the identity on elementary
tensors induces a contraction
\[
\mathbb B(\ell_2)\tmax B\longrightarrow(\mathbb B(\ell_2)\tmin D)/(\mathbb B(\ell_2)\tmin I).
\]
Composing the two contractions gives precisely the criterion for local
$1$-liftability through $q$ in \cite[Proposition~7.48]{Pisier}.

Finally, assume \textup{(iv)} and fix a u.c.p. map $u:A\to\mathbb B(\ell_2)$.
Since $\mathbb B(\ell_2)$ has WEP, $u$ is a WEP map. Choose a free group $\mathbb F$ and
a unital quotient map $q:C^*(\mathbb F)\twoheadrightarrow\mathbb B(\ell_2)$. Given
$z\in\mathbb B(\ell_2)\odot A$, choose a finite-dimensional $E\subset A$ containing
its coefficients and a complete contraction $s_E:E\to C^*(\mathbb F)$
lifting $u|_E$. The full free group algebra has LLP, so
$\mathbb B(\ell_2)\tmax C^*(\mathbb F)=\mathbb B(\ell_2)\tmin C^*(\mathbb F)$. Applying
$\id_{\mathbb B(\ell_2)}\otimes q$ to $(\id_{\mathbb B(\ell_2)}\otimes s_E)z$ gives
\[
\|(\id_{\mathbb B(\ell_2)}\otimes u)z\|_{\mathbb B(\ell_2)\tmax\mathbb B(\ell_2)}
\le \|(\id_{\mathbb B(\ell_2)}\otimes s_E)z\|_{\mathbb B(\ell_2)\tmax C^*(\mathbb F)}
= \|(\id_{\mathbb B(\ell_2)}\otimes s_E)z\|_{\min}
\le \|z\|_{\min}.
\]
This proves \textup{(ii$'$)}.
\end{proof}

\begin{proof}[Proof of Theorem~\ref{thm:QWEP-Pisier}]
By \cite[Theorem~9.67]{Pisier}, \textup{(i)} implies \textup{(ii)},
since every u.c.p. map is a decomposable contraction. Taking
$L_1=L_2=C^*(\F_\infty)$ gives \textup{(ii$'$)}.

Assume \textup{(ii$'$)}. Choose an infinite set $I$ and a unital quotient
map $q:C^*(\F_I)\twoheadrightarrow Q$, where $\F_I$ is free on $I$.
For each $z\in C^*(\F_\infty)\odot C^*(\F_I)$, approximation of its right-hand
coefficients by group polynomials gives a countably infinite $J\subset I$
such that $z\in C^*(\F_\infty)\odot C^*(\F_J)$. Since
$C^*(\F_J)\cong C^*(\F_\infty)$, we may apply \textup{(ii$'$)} to the u.c.p.
map $q|_{C^*(\F_J)}$. By injectivity of the minimal tensor product this gives
\[
\|(\id_{C^*(\F_\infty)}\otimes q)z\|_{\max}\le \|z\|_{\min}.
\]
Thus $q$ is a WEP map, which proves \textup{(iii)}. The implication
\textup{(iii)}$\Rightarrow$\textup{(i)} is \cite[Theorem~9.67]{Pisier}.

Assume \textup{(i)} and write $Q=W/I$ with $W$ having WEP. Let
$p:W\to Q$ be the quotient map. Fix $X\subset Y$ and $u$ as in
\textup{(iv)}. For each finite-dimensional subspace $E\subset X$,
the LLP condition gives a complete contraction $v:E\to W$
with $pv=u|_E$. Indeed, the maximal tensor product preserves quotients,
so the identity on elementary tensors induces a contractive map
\[
\mathbb B(\ell_2)\tmax Q\longrightarrow
(\mathbb B(\ell_2)\tmin W)/(\mathbb B(\ell_2)\tmin I).
\]
Composing with $\id_{\mathbb B(\ell_2)}\otimes u$ gives the local lifting
criterion of \cite[Proposition~7.48]{Pisier}. By WEP,
$i_Wv$ extends to a complete contraction $V:Y\to W^{**}$, where
$i_W:W\to W^{**}$ is the canonical embedding
\cite[Corollary~9.24]{Pisier}. Thus $T_E=p^{**}V:Y\to Q^{**}$ is a
complete contraction agreeing with $i_Q u$ on $E$.

The unit ball of $\mathrm{CB}(Y,Q^{**})$ is compact for pointwise
weak$^*$ convergence. Indeed, it is a closed subset of the compact
product $\prod_{y\in Y}\|y\|B_{Q^{**}}$, since linearity and the
contraction inequalities at every matrix level pass to pointwise
weak$^*$ limits. Order the finite-dimensional subspaces $E\subset X$
by inclusion and take a convergent subnet of $(T_E)$ with limit
$\widetilde u$ in this ball. For each $x\in X$, the net $T_E(x)$ equals
$i_Q u(x)$ as soon as $E$ contains $x$. Thus $\widetilde u|_X=i_Q u$.

Conversely, assume \textup{(iv)}. Choose a free group $\mathbb F$ and a
unital quotient map $q:D=C^*(\mathbb F)\twoheadrightarrow Q$.
Since $D$ has LLP, $\mathbb B(\ell_2)\tmin D=\mathbb B(\ell_2)\tmax D$, so $q$ is an
LLP map. Represent $D$ faithfully and unitally in $\mathbb B(K)$.
\textup{(iv)} gives a complete contraction $T:\mathbb B(K)\to Q^{**}$ with
$T|_D=i_Q q$. In particular $T$ is unital and therefore
completely positive. The factorization $D\hookrightarrow \mathbb B(K)\xrightarrow{T}Q^{**}$
of $i_Q q$ makes $q$ a WEP map by \cite[Theorem~9.57]{Pisier}, proving
\textup{(iii)}.
\end{proof}
We finish this section by proving permanence properties of local embeddability into $C^*(\F_\infty)$.
\begin{proposition}[Permanence of $d_f=1$]\label{prop:df-permanence}
The condition $d_f=1$ has the following permanence properties.
\begin{enumerate}
\item If $d_f(A)=1$ and $B$ is an exact $C^*$-algebra, then
\[
d_f(A\tmin B)=1.
\]
\item If an exact discrete group $G$ acts on a $C^*$-algebra $A$ by
automorphisms $\alpha$, then
\[
d_f(A\rtimes_{\alpha,r}G)=1
\quad\Longleftrightarrow\quad d_f(A)=1.
\]
%Here exactness of $G$ means exactness of $C_r^*(G)$.
\item If $(A_i,\phi_{ji})$ is an inductive system with injective connecting
$*$-homomorphisms and $d_f(A_i)=1$ for every $i$, then
\[
d_f\!\left(\varinjlim_i A_i\right)=1.
\]
\end{enumerate}
\end{proposition}

\begin{proof}
We first record two observations. The condition $d_f=1$ passes to operator
subspaces and is preserved by matrix amplification. For the latter assertion,
a finite-dimensional subspace of $M_k(A)$ lies in $M_k(E)$ for some
finite-dimensional $E\subset A$; amplify an almost completely isometric
embedding of $E$ into $C^*(\F_\infty)$ and use $d_f(M_k(C^*(\F_\infty)))=1$.
This last equality follows from the LLP of $M_k(C^*(\F_\infty))$.

The condition also passes through completely contractive approximate
factorizations of the identity. More precisely, suppose
\[
X\xrightarrow{\varphi_\lambda}Y_\lambda
\xrightarrow{\psi_\lambda}X,
\qquad
\|\varphi_\lambda\|_{\cb},\|\psi_\lambda\|_{\cb}\le1,
\qquad d_f(Y_\lambda)=1,
\]
and $\psi_\lambda\varphi_\lambda\to\id_X$ pointwise in norm.
For every operator space $F$ and $z\in X\odot F$, functoriality gives
\[
\|((\psi_\lambda\varphi_\lambda)\otimes\id_F)z\|_M
\le \|(\varphi_\lambda\otimes\id_F)z\|_{\min}
\le \|z\|_{\min}.
\]
The left-hand tensors converge to $z$ in the $M$-norm, so
\cite[Theorem~4.5]{JP} implies $d_f(X)=1$.

\smallskip\noindent
\textup{(i)} Let $E\subset A\odot B$ be finite dimensional, and choose a
finite-dimensional coefficient space $F\subset B$ with $E\subset A\odot F$.
Exactness of $B$ implies that, for every $\delta>0$, there are $k$ and an
isomorphism $v:F\to v(F)\subset M_k$ satisfying
$\|v\|_{\cb}\|v^{-1}\|_{\cb}<1+\delta$.
Minimal tensorization gives an embedding of $E$ into $M_k(A)$ with the
same distortion bound. Since $d_f(M_k(A))=1$, we obtain
$d_f(E)\le1+\delta$, and hence $d_f(E)=1$.
An arbitrary finite-dimensional subspace of $A\tmin B$ is an arbitrarily
small cb perturbation of one in $A\odot B$: approximate a basis, using the
estimate in part~\textup{(iii)} below. This proves~\textup{(i)}.

\smallskip\noindent
\textup{(ii)} Assume $d_f(A)=1$. By the amenable-action characterization of
exact discrete groups, there is a nonempty compact $G$-space $X$ on which
$G$ acts topologically amenably; see \cite[Chapter~5]{BO}.
Put $D=A\tmin C(X)$ and let $\beta$ be the diagonal action. Part~\textup{(i)}
gives $d_f(D)=1$.

We spell out the approximate factorizations supplied by amenability.
Write $\gamma$ for the action on $C(X)$. There is a net of finitely
supported families $h=(h_t)_{t\in G}$ of nonnegative functions in $C(X)$
such that
\[
\sum_t h_t^2=1,
\qquad
c_h(g):=\sum_t h_{gt}\gamma_g(h_t)\longrightarrow1
\quad\text{in norm for every }g\in G.
\]
Identify these functions with central multipliers of $D$. Let
$F=\operatorname{supp}(h)$ and denote the canonical unitaries in the
multiplier algebra of $D\rtimes_{\beta,r}G$ by $w_g$. Regular compression
and row compression, respectively, define completely positive contractions
\[
\Phi_h:D\rtimes_{\beta,r}G\longrightarrow M_F(D),
\qquad
\Psi_h:M_F(D)\longrightarrow D\rtimes_{\beta,r}G,
\]
\[
\begin{aligned}
\Phi_h(dw_g)_{s,t}&=\delta_{s,gt}\,\beta_s^{-1}(d),\\
\Psi_h([d_{s,t}])&=\sum_{s,t\in F}h_sw_s d_{s,t}w_t^*h_t.
\end{aligned}
\]
Indeed, the row $(h_sw_s)_{s\in F}$ is contractive because
$\sum_s h_sw_sw_s^*h_s=1$. Direct calculation gives
\[
\Psi_h\Phi_h(dw_g)=c_h(g)dw_g.
\]
Thus $\Psi_h\Phi_h\to\id$ pointwise in norm. Since $d_f(M_F(D))=1$,
the preceding observation yields $d_f(D\rtimes_{\beta,r}G)=1$.
The equivariant embedding $a\mapsto a\otimes1$ induces an injective
$*$-homomorphism
\[
A\rtimes_{\alpha,r}G\hookrightarrow D\rtimes_{\beta,r}G,
\]
so $d_f(A\rtimes_{\alpha,r}G)=1$. Conversely, $A$ is a $C^*$-subalgebra
of its reduced crossed product, which proves the reverse implication.

\smallskip\noindent
\textup{(iii)} Identify the stages with subalgebras of
$A=\overline{\bigcup_i A_i}$. Fix a finite-dimensional $E\subset A$ and a
biorthogonal basis $(x_j,f_j)_{j=1}^d$ of $E$. Given $0<\delta<1$, choose
elements $y_j$ in a common stage $A_i$ such that
$\sum_j\|f_j\|\,\|y_j-x_j\|<\delta$. The map
$T:E\to A_i$, $T(x)=\sum_jf_j(x)y_j$, satisfies
\[
\|T-\iota_E\|_{\cb}
\le\sum_j\|f_j\|\,\|y_j-x_j\|<\delta,
\]
because a scalar-valued functional has cb norm equal to its norm.
Consequently $T$ is injective and
\[
d_f(E)\le
d_{\cb}(E,T(E))\,d_f(T(E))
\le\frac{1+\delta}{1-\delta}.
\]
Letting $\delta\downarrow0$ and then taking the supremum over $E$ proves
\textup{(iii)}.
\end{proof}

\begin{remark}[Failure for non-injective inductive limits]
\label{rem:df-noninjective-limits}
The injectivity assumption in Proposition~\ref{prop:df-permanence}\textup{(iii)}
cannot be omitted, even when all building blocks have LP and all connecting
maps are surjective. To see this, choose a separable $C^*$-algebra $B$ with
$d_f(B)=\infty$. Such a $B$ exists: by \cite{JP}, there is a
finite-dimensional operator space $E\subset \mathbb B(\ell_2)$ which
does not locally embed into $C^*(\F_\infty)$. Taking $B=C^*(E,1)$,
Proposition~\ref{prop:df-dichotomy} gives $d_f(B)=\infty$. Its cone
\[
CB=C_0((0,1])\tmin B
\]
contains a completely isometric operator-space copy of $B$, through
$b\mapsto[t\mapsto tb]$. Hence $d_f(CB)=d_f(B)=\infty$.
On the other hand, the cone theorem of Loring--Shulman, recalled in
\cite[Theorem~4.1]{Thiel}, writes $CB$ as an inductive limit of separable
projective $C^*$-algebras $P_n$ with surjective connecting maps; see also the proof of
\cite[Proposition~5.7]{ES}. Each $P_n$ has LP, and therefore $d_f(P_n)=1$.
Thus
\[
d_f(P_n)=1\quad(n\ge1),
\qquad
d_f\!\left(\varinjlim_n P_n\right)=\infty.
\]
Unitizing this system gives the same counterexample with unital building
blocks and unital surjective connecting maps.
\end{remark}

\section{Ambient CP-stability}
\label{sec:ACP}

CP-stability was introduced by Goldbring and Sinclair
\cite[Definition~2.38]{GS} and further studied by Sinclair \cite{Sinclair}.
We consider an ambient variant, called \emph{ambient CP-stability}
(\emph{ACP-stability} for short), which characterizes $d_f(A)=1$.
The test maps are defined on finite-dimensional operator systems in an
ambient $\mathbb B(H)$.

An operator system is a unital self-adjoint linear subspace of a unital
$C^*$-algebra, with its inherited matrix order. For a linear map $\varphi$
between operator systems, write
\[
\|\varphi\|_k=\|\id_{M_k}\otimes\varphi\|.
\]
\begin{definition}[Ambient CP-stability]\label{def:ACP}
A unital $C^*$-algebra $A$ is \emph{ACP-stable} if, for every faithful
unital representation $j:A\hookrightarrow \mathbb B(H)$, every finite-dimensional
operator system $E\subset A$, and every $\varepsilon>0$, there exist a
finite-dimensional operator system
\[
j(E)\subset F\subset \mathbb B(H),
\qquad k\in\N,\quad\delta>0,
\]
such that the following holds for every unital $C^*$-algebra $B$:
every unital linear map $\varphi:F\to B$ satisfying
$\|\varphi\|_k<1+\delta$ admits a u.c.p. map $\psi:E\to B$ with
\begin{equation}\label{eq:ACP}
\|\psi-(\varphi\circ j)|_E\|_{\cb}<\varepsilon.
\end{equation}
The data $F,k,\delta$ may depend on $j,E,\varepsilon$, but not on
$B$ or $\varphi$. In particular, $F$ is not required to be contained in
$j(A)$.
\end{definition}

The cb-norm in \eqref{eq:ACP} can equivalently be replaced by the ordinary
operator norm. Indeed, for a fixed finite-dimensional operator space $E$
there is a constant $c_E$ such that
\begin{equation}\label{eq:finite-dimensional-cb}
\|v\|_{\cb}\le c_E\|v\|
\end{equation}
for every linear map $v:E\to Y$ and every operator space $Y$.
To see this, choose a basis $(e_r)_{r=1}^d$ of $E$ and its coordinate
functionals $(f_r)_{r=1}^d$. The decomposition
$v=\sum_r f_r(\,\cdot\,)v(e_r)$ gives
$\|v\|_{\cb}\le\sum_r\|f_r\|\|v(e_r)\|$.
Thus one may take $c_E=\sum_r\|f_r\|\|e_r\|$.

It suffices to test Definition~\ref{def:ACP} in one faithful unital
representation. Indeed, given two such representations
$j_r:A\hookrightarrow \mathbb B(H_r)$, $r=1,2$, Arveson's extension theorem gives
a u.c.p. map
\[
\Theta:\mathbb B(H_1)\longrightarrow \mathbb B(H_2),\qquad
\Theta\circ j_1=j_2.
\]
If the condition in Definition~\ref{def:ACP} holds for $j_1$, choose its
control system $F_1$ and put $F_2=\Theta(F_1)$. This is a finite-dimensional
operator system containing $j_2(E)$. For a test map
$\varphi_2:F_2\to B$, the map
$\varphi_1=\varphi_2\circ\Theta|_{F_1}$ is unital and satisfies
$\|\varphi_1\|_k\le\|\varphi_2\|_k$, while
$\varphi_1\circ j_1|_E=\varphi_2\circ j_2|_E$.
Thus the same $k,\delta$ and the resulting correction $\psi:E\to B$
work for $j_2$.

\Needspace{6\baselineskip}
\begin{theorem}\label{thm:ACP}
A unital $C^*$-algebra $A$ is ACP-stable if and only if $d_f(A)=1$.
\end{theorem}

\begin{proof}
Fix a faithful unital representation $j:A\hookrightarrow \mathbb B(H)$ and
identify $A$ with $j(A)$.

Assume $d_f(A)=1$. If $A$ fails to be ACP-stable, there are a
finite-dimensional operator system $E\subset A$ and $\varepsilon_0>0$
for which Definition~\ref{def:ACP} fails.
Let $\mathcal I$ be the directed set of pairs $\alpha=(F,n)$, where
$E\subset F\subset \mathbb B(H)$ is a finite-dimensional operator system and
$n\ge1$, ordered by inclusion in the first coordinate and the usual order
in the second. For each $\alpha=(F,n)$ choose a unital $C^*$-algebra
$B_\alpha$ and a unital linear map $\varphi_\alpha:F\to B_\alpha$ such
that
\begin{equation}\label{eq:ACP-counterexamples}
\|\varphi_\alpha\|_n<1+1/n,
\qquad
\inf_{\psi\in\mathrm{UCP}(E,B_\alpha)}
\|\psi-\varphi_\alpha|_E\|_{\cb}\ge\varepsilon_0.
\end{equation}
Choose an ultrafilter $\omega$ on $\mathcal I$ containing every tail, and
form the norm ultraproduct
\[
Q=\left(\prod_{\alpha\in\mathcal I}B_\alpha\right)/I_\omega,
\qquad
I_\omega=\{(b_\alpha):\lim_\omega\|b_\alpha\|=0\}.
\]
For $x\in \mathbb B(H)$ define a bounded family by
\[
b_\alpha(x)=
\begin{cases}
\varphi_\alpha(x),&x\in F,\quad\alpha=(F,n),\\
0,&x\notin F,
\end{cases}
\qquad
\Phi(x)=[(b_\alpha(x))]_\omega.
\]
The coordinate assignments need not be linear on all of $\mathbb B(H)$.
Nevertheless, their classes define a unital linear map $\Phi:\mathbb B(H)\to Q$:
any finite set of elements is contained in $F$ on an $\omega$-large set.
At each fixed matrix level $r$, on an $\omega$-large set the entries of a
given matrix lie in $F$ and $n\ge r$. The first inequality in
\eqref{eq:ACP-counterexamples} then gives $\|\Phi\|_r\le1$.
Thus $\Phi$ is unital and completely contractive, hence u.c.p.

The restriction $u=\Phi|_A$ is a WEP map. Indeed, $\id_{C^*(\F_\infty)}\otimes u$
factors contractively as
\[
C^*(\F_\infty)\tmin A\longrightarrow
C^*(\F_\infty)\tmin \mathbb B(H)=C^*(\F_\infty)\tmax \mathbb B(H)
\xrightarrow{\ \id_{C^*(\F_\infty)}\otimes\Phi\ }C^*(\F_\infty)\tmax Q,
\]
where the equality follows from WEP of $\mathbb B(H)$. By
Theorem~\ref{thm:df-WEP-maps}\textup{(iv)}, $u$ is locally $1$-liftable.
Since $u$ is unital and completely positive, \cite[Corollary~9.47]{Pisier}
gives a u.c.p. lifting
\[
\Psi:E\longrightarrow\prod_{\alpha\in\mathcal I}B_\alpha.
\]
Write $\psi_\alpha$ for its coordinate maps. Each $\psi_\alpha$ is
u.c.p., and
\[
\lim_\omega\|\psi_\alpha(e)-\varphi_\alpha(e)\|=0
\qquad(e\in E).
\]
The basis estimate in the proof of \eqref{eq:finite-dimensional-cb}
then gives $\lim_\omega\|\psi_\alpha-\varphi_\alpha|_E\|_{\cb}=0$.
This contradicts \eqref{eq:ACP-counterexamples}.

Conversely, assume that $A$ is ACP-stable. Choose a unital
surjection $q:C^*(\F_\Lambda)\twoheadrightarrow \mathbb B(H)$ from a full free
group $C^*$-algebra. Fix a finite-dimensional operator system
$E\subset A$ and $0<\varepsilon<1$. Choose $F,k,\delta$ as in
Definition~\ref{def:ACP}. By the finite-level lifting theorem
\cite{RS}, applied as in \cite[proof of Proposition~2.5]{Sinclair},
the inclusion $F\hookrightarrow \mathbb B(H)$ admits a unital $2k$-positive
lifting $\varphi:F\to C^*(\F_\Lambda)$ through $q$. Applying
$2k$-positivity to $\left(\begin{smallmatrix}1&x\\x^*&1\end{smallmatrix}\right)$
for $x$ in the unit ball of $M_k(F)$ gives $\|\varphi\|_k\le1$.
Thus ACP-stability provides a u.c.p. map
$s:E\to C^*(\F_\Lambda)$ with $\|s-\varphi|_E\|_{\cb}<\varepsilon$.

Since $q\varphi|_E$ is the inclusion of $E$ into $\mathbb B(H)$, for every
$r\ge1$ and $x\in M_r(E)$ we have
\[
(1-\varepsilon)\|x\|
\le \|(\id_{M_r}\otimes qs)x\|
\le \|(\id_{M_r}\otimes s)x\|
\le \|x\|.
\]
Approximating a basis of $s(E)$ by group polynomials shows that $s(E)$
lies in the subalgebra generated by countably many free generators.
This subalgebra embeds into $C^*(\F_\infty)$, so
$d_f(E)\le (1-\varepsilon)^{-1}$. Letting $\varepsilon\downarrow0$
gives $d_f(E)=1$. Every finite-dimensional subspace of $A$ lies in a
finite-dimensional operator system, hence $d_f(A)=1$.
\end{proof}

\begin{remark}[Comparison with the LLP]\label{rem:ACP-LLP}
If one requires $F\subset j(A)$ in Definition~\ref{def:ACP}, one obtains
ordinary CP-stability, which is equivalent to the LLP
\cite[Proposition 2.5]{Sinclair}; the use of the cb-norm in the conclusion
does not change the property by \eqref{eq:finite-dimensional-cb}.
Allowing $F\subset \mathbb B(H)$ therefore gives ACP-stability whenever $A$
has the LLP.
\end{remark}

\section{Amenable traces and relative local lifting}\label{sec:tracial-lifting}

%Local lifting relative to an amenable trace gives a tracial counterpart
%of the implication
It is well known that if a unital $C^*$-algebra $A$ has WEP and locally embedds into $C^*(\F_\infty)$, then it
has LLP.
This principle remains if one replaces the WEP by the presence of an amenable trace, but the conclusion should be replaced by a local lifting property relative to the trace, as we will see in Proposition~\ref{prop:trace-LLP} below.
This applies to any amenable trace and we will use it
to prove that $d_f(C^*(\SL_3(\Z)))=\infty$.

\begin{definition}\label{def:trace-LLP}
Let $A$ be a unital $C^*$-algebra and let $\tau$ be a tracial state on $A$.
We say that $A$ has $\tau\text{-}\mathrm{LLP}$ if the following holds. For every
finite von Neumann algebra $N$ with a faithful normal tracial state
$\sigma$, every unital $*$-homomorphism
\[
\pi:A\longrightarrow N,\qquad \sigma\circ\pi=\tau,
\]
every unital surjective $*$-homomorphism $q:B\twoheadrightarrow N$
from a unital $C^*$-algebra $B$, and every finite-dimensional operator
system $E\subset A$, there exists a u.c.p. map
\[
s:E\longrightarrow B,\qquad q\circ s=\pi|_E.
\]
The notation $\tau\text{-}\mathrm{LLP}$ is used here for this property of the
pair $(A,\tau)$; it does not denote the usual LLP of $A$.
\end{definition}

The usual LLP implies $\tau\text{-}\mathrm{LLP}$ for every tracial state $\tau$.
In Definition~\ref{def:trace-LLP}, however, only representations realizing
the specified trace are tested. Neither $\pi$ nor $\tau$ is assumed
faithful, and no normality is imposed on $q$ or on the liftings.

\begin{proposition}\label{prop:trace-LLP}
Let $A$ be a unital $C^*$-algebra and let $\tau$ be an amenable tracial
state on $A$. If $d_f(A)=1$, then $A$ has $\tau\text{-}\mathrm{LLP}$.
\end{proposition}

\begin{proof}
Fix $N,\sigma,\pi,q,E$ as in Definition~\ref{def:trace-LLP}, and choose
a faithful unital representation $j:A\hookrightarrow \mathbb B(H)$.
By the amenable-trace extension theorem
\cite[Theorem 3.1.6(5)]{Brown}, there exists a u.c.p. map
\begin{equation}\label{eq:trace-weak-expectation}
\Psi_\tau:\mathbb B(H)\longrightarrow M_\tau:=\pi_\tau(A)'',
\qquad \Psi_\tau(j(a))=\pi_\tau(a)\quad(a\in A).
\end{equation}
Set $P=\pi(A)''\subset N$. The equality $\sigma\pi=\tau$, together
with the faithfulness and normality of $\sigma|_P$, identifies the GNS
von Neumann algebra $M_\tau$ with $P$. More precisely, uniqueness of the
GNS representation gives a trace-preserving normal $*$-isomorphism
\[
\theta:M_\tau\longrightarrow P,
\qquad \theta(\pi_\tau(a))=\pi(a).
\]
Consequently, $\widetilde\pi=\theta\Psi_\tau$, viewed as a map into
$N$, is u.c.p. and
\begin{equation}\label{eq:trace-factorization}
\pi=\widetilde\pi\circ j:
A\xrightarrow{\ j\ }\mathbb B(H)
\xrightarrow{\ \widetilde\pi\ }N=B/\ker q.
\end{equation}
This factorization through $\mathbb B(H)$ makes $\pi$ a WEP map, since
$C^*(\F_\infty)\tmin \mathbb B(H)=C^*(\F_\infty)\tmax \mathbb B(H)$ and $\widetilde\pi$ is u.c.p.
As $d_f(A)=1$, Theorem~\ref{thm:df-WEP-maps}\textup{(iv)} gives local
$1$-liftability through $q$. By \cite[Corollary~9.47]{Pisier}, $\pi|_E$
therefore has a u.c.p. lifting $s:E\to B$.
\end{proof}

The analogy with WEP is at the level of the extension step. For a
faithful $j:A\hookrightarrow \mathbb B(H)$, WEP provides a u.c.p. map
$\mathbb B(H)\to A^{**}$ extending the canonical inclusion of $A$ into its
bidual. Amenability of $\tau$ instead provides
\eqref{eq:trace-weak-expectation}, whose target is the particular
von Neumann quotient $M_\tau$ of $A^{**}$. This suffices for
\eqref{eq:trace-factorization} and hence for the liftings tested by
$\tau\text{-}\mathrm{LLP}$. Proposition~\ref{prop:trace-LLP} does not assert
that an algebra with $d_f(A)=1$ and one amenable trace has the usual LLP.

\begin{theorem}\label{thm:SL3-tracial-lifting}
Let $\Gamma$ be either $\Z^2\rtimes\SL_2(\Z)$ or $\SL_3(\Z)$.
Then the full group $C^*$-algebra $C^*(\Gamma)$ does not locally embed
into $C^*(\F_\infty)$, i.e., $d_f(C^*(\Gamma))=\infty$.
\end{theorem}

\begin{proof}
In this proof, we follow \cite[ISW20]{ISW} but we incorporate our amenable trace trick. First take $\Gamma=\Z^2\rtimes\SL_2(\Z)$, and put $\Lambda=\Z^2$
and $A=C^*(\Gamma)$. We first describe the
representation to which Proposition~\ref{prop:trace-LLP} will apply.
For $a,b\in\Z^2$, write $\kappa(a,b)=a_1b_2-a_2b_1$, and set
\[
c((a,S),(b,T))=\kappa(a,Sb)
\qquad ((a,S),(b,T)\in\Gamma).
\]
The $\SL_2(\Z)$-invariance of $\kappa$ implies that $c$ is a normalized
integer-valued $2$-cocycle. Choose strictly increasing odd primes
$p_n$ and let
\[
c_n(g,h)=\exp\bigl(2\pi i c(g,h)/p_n\bigr).
\]
These multipliers converge pointwise to $1$. Their restrictions to
$\Lambda$ are not coboundaries: a scalar coboundary on an abelian group
is symmetric, whereas $c_n(e_1,e_2)=\exp(2\pi i/p_n)$ and
$c_n(e_2,e_1)=\exp(-2\pi i/p_n)$, where $e_1,e_2$ are the standard
basis vectors of $\Lambda=\Z^2$.

Let $r_n:\Gamma\twoheadrightarrow Q_n$ be reduction modulo $p_n$, where
\[
Q_n=(\Z/p_n\Z)^2\rtimes\SL_2(\Z/p_n\Z),
\qquad d_n=|Q_n|.
\]
The cocycle $c_n$ descends to a cocycle $\mu_n$ on $Q_n$. Its twisted
left regular representation, pulled back to $\Gamma$, is
\[
v_n(g)\delta_x=\mu_n(r_n(g),x)\delta_{r_n(g)x}
\qquad(x\in Q_n).
\]
Thus $v_n(g)\in\U(M_{d_n})$ and
\begin{equation}\label{eq:trace-projective-reps}
v_n(g)v_n(h)=c_n(g,h)v_n(gh).
\end{equation}
This is the finite-quotient construction of \cite[Section 3.2]{ISW}.
For every fixed $g\ne e$, one has $r_n(g)\ne e$ for all sufficiently
large $n$. The corresponding permutation has no fixed points, and hence
\begin{equation}\label{eq:trace-projective-limit}
\lim_n\tr_{d_n}(v_n(g))=\delta_{g,e}\qquad(g\in\Gamma).
\end{equation}

The \emph{untwisted} left regular representations of $Q_n$, composed
with $r_n$, are genuine finite-dimensional representations of $\Gamma$.
Their normalized traces also converge pointwise to the canonical trace
$\tau_\Gamma(u_g)=\delta_{g,e}$. It follows that $\tau_\Gamma$ is an
amenable trace on the full algebra $A$; see
\cite[Theorem 3.1.6]{Brown}. Pointwise convergence on the canonical
unitaries extends to all of $A$ by density and the norm-one bound on
these states.

Fix a free ultrafilter $\omega$ and form the tracial ultraproduct
\[
N=\Bigl(\prod_nM_{d_n}\Bigr)/I_\omega,
\qquad
I_\omega=\bigl\{(x_n):\lim_{n\to\omega}
             \|x_n\|_{2,\tr_{d_n}}=0\bigr\}.
\]
Write $q_\omega:\prod_nM_{d_n}\twoheadrightarrow N$ for the quotient
map and $\tau_\omega$ for the faithful normal tracial state on $N$.
Equation~\eqref{eq:trace-projective-reps} and $c_n\to1$ show that
\[
\pi:A\longrightarrow N,
\qquad \pi(u_g)=q_\omega((v_n(g))_n),
\]
is a unital $*$-homomorphism. By
\eqref{eq:trace-projective-limit}, it satisfies
\begin{equation}\label{eq:trace-obstructive-representation}
\tau_\omega\circ\pi=\tau_\Gamma.
\end{equation}

Assume, towards a contradiction, that $d_f(A)=1$.
Proposition~\ref{prop:trace-LLP} implies that $A$ has
$\tau_\Gamma\text{-}\mathrm{LLP}$. We show that even one local lifting of
$\pi$ through $q_\omega$ is impossible.

\Needspace{9\baselineskip}
The relative property $(T)$ of $(\Gamma,\Lambda)$ supplies a finite
set $F\subset\Gamma$ and $\delta>0$ such that every projective unitary
representation $w$ whose multiplier restricts to a non-coboundary on
$\Lambda$ satisfies
\begin{equation}\label{eq:trace-projective-rigidity}
\max_{g\in F}\|w(g)\xi-\xi\|\ge\delta\|\xi\|
\quad\text{for every vector $\xi$}.
\end{equation}
This estimate follows from \cite[Theorem~2.3]{ISW}, as in the proof
of Claim~3.1 there.
Set
\[
E=\spanop\{1,u_g,u_g^*:g\in F\}\subset A.
\]
A u.c.p. lifting $s:E\to\prod_nM_{d_n}$ of $\pi|_E$ has coordinate
maps $s_n:E\to M_{d_n}$ satisfying
\begin{equation}\label{eq:trace-local-approximation}
\lim_{n\to\omega}
\|s_n(u_g)-v_n(g)\|_{2,\tr_{d_n}}=0
\qquad(g\in F).
\end{equation}
Arveson's extension theorem extends each $s_n$ to a u.c.p. map
$\varphi_n:A\to M_{d_n}$. Choose a Stinespring representation
\[
\varphi_n(a)=V_n^*\rho_n(a)V_n,
\qquad V_n:\mathbb C^{d_n}\longrightarrow K_n\text{ an isometry},
\]
where $\rho_n:A\to \mathbb B(K_n)$ is a unital $*$-representation.

On the Hilbert space $\mathcal S_2(\mathbb C^{d_n},K_n)$ of Hilbert--Schmidt
operators, with its usual unnormalized Hilbert--Schmidt norm, define
\[
w_n(g)T=\rho_n(u_g)T v_n(g)^*,
\qquad \xi_n=d_n^{-1/2}V_n.
\]
Each $w_n$ is a projective unitary representation with multiplier
$\overline{c_n}$, and $\|\xi_n\|=1$. A direct calculation gives
\begin{align}
\|w_n(g)\xi_n-\xi_n\|^2
&=2-2\Ree\tr_{d_n}\bigl(\varphi_n(u_g)v_n(g)^*\bigr)\notag\\
&\le 2\|\varphi_n(u_g)-v_n(g)\|_{2,\tr_{d_n}}.
\label{eq:trace-Stinespring-defect}
\end{align}
For $g\in F$, the last expression tends to zero along $\omega$ by
\eqref{eq:trace-local-approximation}. Since $F$ is finite,
\[
\lim_{n\to\omega}\max_{g\in F}
\|w_n(g)\xi_n-\xi_n\|=0.
\]
But $\overline{c_n}|_\Lambda$ is not a coboundary for any $n$, so
\eqref{eq:trace-projective-rigidity} makes every term at least $\delta$.
This contradiction rules out $d_f(C^*(\Gamma))=1$.
Proposition~\ref{prop:df-dichotomy} now gives $d_f(C^*(\Gamma))=\infty$.

Finally, the embedding
\[
\Gamma\hookrightarrow\SL_3(\Z),
\qquad (a,S)\longmapsto
\begin{pmatrix}S&a\\0&1\end{pmatrix},
\]
induces a faithful $*$-homomorphism of the full group $C^*$-algebras.
Monotonicity of $d_f$ under complete isometries therefore gives
\[
\infty=d_f(C^*(\Gamma))\le d_f(C^*(\SL_3(\Z))).
\]
\end{proof}

\begin{remark}\label{rem:trace-versus-LLP-obstruction}
For $\Gamma=\Z^2\rtimes\SL_2(\Z)$ the proof establishes the failure of
$\tau_\Gamma\text{-}\mathrm{LLP}$ for $C^*(\Gamma)$. The conclusion
$d_f(C^*(\Gamma))=\infty$ follows from this failure, the amenability
of $\tau_\Gamma$, and Proposition~\ref{prop:df-dichotomy}.
The extra input that permits the use of $d_f(A)=1$ in the
Ioana--Spaas--Wiersma obstruction is precisely the ambient extension
of the representation in \eqref{eq:trace-obstructive-representation}.
\end{remark}

\begin{remark}\label{rem:SL3-polynomial-rings}
For every prime $p\ge3$, the same argument gives
\[
d_f\bigl(C^*(\SL_3(\F_p[t]))\bigr)=\infty.
\]
Indeed, put $R=\F_p[t]$ and $\Gamma=R^2\rtimes\SL_2(R)$.
As in \cite[Corollary~B and Section~3.2]{ISW}, use the determinant
cocycle and the additive characters of $R/(t^{n+1})$ obtained by
extracting the coefficient of $t^n$ and identifying $\F_p$ with the
$p$th roots of unity. The resulting scalar cocycles converge pointwise
to $1$ and are not symmetric on $R^2$ when $p\ge3$, so their
restrictions are not coboundaries. Their twisted regular
representations on the finite congruence quotients have normalized
traces converging to the canonical trace $\tau_\Gamma$.
The untwisted regular representations of these quotients show that
$\tau_\Gamma$ is amenable. The relative property $(T)$ of
$(\Gamma,R^2)$, recalled in \cite[Section~3.2]{ISW}, and
Proposition~\ref{prop:trace-LLP} therefore allow us to repeat the proof
of Theorem~\ref{thm:SL3-tracial-lifting}, giving
$d_f(C^*(\Gamma))=\infty$. The block embedding
$\Gamma\hookrightarrow\SL_3(R)$ yields the assertion.
This argument does not cover $p=2$.
\end{remark}

%\Needspace{12\baselineskip}
\begin{remark}[The Junge--Pisier construction and quotients of the full group algebra]
\label{rem:JP-quotients}
The use of $\SL_3(\Z)$ in this context already appears in
\cite[Remark~2.10 and the subsequent paragraph, arXiv version]{JP}.
Junge and Pisier's argument builds on Voiculescu's construction in
\cite{Voiculescu}, more precisely on the proof of Proposition~3 there.
It is important, however, to distinguish the resulting operator spaces
from subspaces of the full group $C^*$-algebra with its original matrix norms.
Let $G=\SL_3(\Z)$ and $A=C^*(G)$. Starting from pairwise disjoint
finite-dimensional representations $(\pi_k)$, Junge and Pisier consider
\[
\pi_\Omega=\bigoplus_{k\in\Omega}\pi_k,\qquad
T_j^\Omega=\pi_\Omega(g_j),\qquad
M_\Omega=\spanop\{T_1^\Omega,\ldots,T_n^\Omega\},
\]
where $g_1,\ldots,g_{n-1}$ generate $G$ and $g_n=e$.
Their argument gives a nonseparable family of such operator spaces for
the metric $\log d_{\cb}$. Since $G$ has two generators, the construction
can be carried out in dimension three.

These spaces lie naturally in the image algebras
\[
A_\Omega=C^*(\pi_\Omega(G))\cong A/\ker\pi_\Omega.
\]
If $E=\spanop\{u_{g_1},\ldots,u_{g_n}\}\subset A$, then
$M_\Omega=\pi_\Omega(E)$, but the restriction of $\pi_\Omega$ to $E$
is only known to be completely contractive. For every $r\ge1$ and
$a_1,\ldots,a_n\in M_r$,
\[
\left\|\sum_j a_j\otimes T_j^\Omega\right\|
\le \left\|\sum_j a_j\otimes u_{g_j}\right\|.
\]
There is no assertion that these inequalities are equalities. In fact,
the variation of the matrix norms on the image spaces is essential to
the construction.

One can nevertheless deduce that $d_f(A_\Omega)=\infty$ for some $\Omega$.
Indeed, by \cite[Proposition~2.6(i)]{JP}, the set of three-dimensional
subspaces of the separable algebra $C^*(\F_\infty)$ is
separable for $\log d_{\cb}$. Its closure is therefore separable as well,
and contains every three-dimensional operator space $F$ with $d_f(F)=1$.
It cannot contain the entire nonseparable family above. Consequently,
some $M_\Omega$ does not locally embed into $C^*(\F_\infty)$, and
neither does $A_\Omega$. Proposition~\ref{prop:df-dichotomy} then gives
\[
d_f(A_\Omega)=\infty.
\]
Thus the Junge--Pisier construction already yields quotients of
$C^*(\SL_3(\Z))$ which fail local representability in $C^*(\F_\infty)$.

This conclusion does not by itself imply $d_f(A)=\infty$: the property
$d_f=1$ does not pass to arbitrary $C^*$-algebra quotients. The obstruction
to such an inference is particularly transparent here. Precomposing
every $\pi_\Omega$ with a fixed surjection $\F_2\twoheadrightarrow G$
produces exactly the same spaces $M_\Omega$ from representations of
$\F_2$, whereas $d_f(C^*(\F_2))=1$. Accordingly, the distinction between
the conclusion of the earlier construction and the result established
here concerns the full algebra $C^*(\SL_3(\Z))$ itself, rather than
certain image algebras. The group and its rigidity already play a role
in \cite{JP}; their use alone is not a new feature.
\end{remark}

We now introduce Thom's second group.
Let $p$ be a prime and put $R=\Z[1/p]$. Consider the group
\[
\widetilde K_p=
\left\{
\begin{pmatrix}
1&*&*&*\\
0&S&*&*\\
0&0&T&*\\
0&0&0&1
\end{pmatrix}
:S,T\in\SL_3(R)
\right\}\subset\SL_8(R),
\]
where the block sizes are $1,3,3,1$ and the entries in the off-diagonal
blocks are arbitrary elements of $R$. Let
\[
K_p=\widetilde K_p/\{I_8+mE_{18}:m\in\Z\},
\]
where $E_{18}$ is the standard matrix unit. This is the group constructed
in \cite[\S3]{Thom}.

\begin{corollary}[Thom's second group]\label{cor:Thom-second}
For every prime $p$, the full group $C^*$-algebra $C^*(K_p)$ does not locally embed into
$C^*(\F_\infty)$, i.e., $d_f(C^*(K_p))=\infty$.
\end{corollary}

\begin{proof}
The homomorphism
\[
\SL_3(\Z)\longrightarrow\widetilde K_p,
\qquad S\longmapsto\operatorname{diag}(1,S,I_3,1),
\]
is injective, and its image meets the central subgroup
$\{I_8+mE_{18}:m\in\Z\}$ trivially. It therefore induces an embedding
$\SL_3(\Z)\hookrightarrow K_p$, and hence a faithful $*$-homomorphism
$C^*(\SL_3(\Z))\hookrightarrow C^*(K_p)$. By
Theorem~\ref{thm:SL3-tracial-lifting} and monotonicity under complete
isometries,
\[
\infty=d_f(C^*(\SL_3(\Z)))\le d_f(C^*(K_p)).
\]
\end{proof}

\section{Thom's first group}\label{sec:Thom-first}
In this last section, we want to illustrate that the amenable trace trick can apply to other groups whose full $C^*$-algebra fails the LLP.
Thom constructed, for each prime $p$, a group $G_p$ with property $(T)$
whose full group $C^*$-algebra does not have the LLP
\cite[Theorem~1.4 and the discussion following it, p.~198]{Thom}.
We recall its construction from \cite[\S2.1]{Thom}.
Fix a prime $p$ and put $R=\F_p[t,t^{-1}]$. Define
\[
\widetilde G_p=
\left\{
\begin{pmatrix}
1&x&z\\
0&S&y\\
0&0&1
\end{pmatrix}
: S\in\SL_3(R),\ x\in M_{1,3}(R),\
y\in M_{3,1}(R),\ z\in R
\right\}.
\]
The block sizes are $1,3,1$. Write $e_{ij}(a)=I_5+aE_{ij}$ for
$i\ne j$, where $E_{ij}$ is the standard matrix unit, and set
\[
Z_p=\{e_{15}(a):a\in R\},\qquad
C_p=\{e_{15}(a):a\in\F_p[t]\},\qquad
G_p=\widetilde G_p/C_p.
\]
Here $Z_p$ is the centre of $\widetilde G_p$, and $G_p$ has property
$(T)$ by \cite[Lemma~2.1]{Thom}.

\begin{theorem}\label{thm:Thom-first}
For every prime $p$, the full group $C^*$-algebra $C^*(G_p)$ does not
locally embed into $C^*(\F_\infty)$, i.e.,
\[
d_f(C^*(G_p))=\infty.
\]
\end{theorem}

\begin{proof}
Put $N_p=Z_p/C_p$. We first construct finite-dimensional projective
representations of $G_p$. For $n\ge1$, let $R_n=R/(t^{2n+1}-1)$.
Every element of this finite ring has a unique representative
$\sum_{j=0}^{2n}a_jt^j$ with $a_j\in\F_p$. Define the additive character
\[
\chi_n\left(\sum_{j=0}^{2n}a_jt^j\right)
=\exp(2\pi i a_n/p).
\]
Let $r_n:\widetilde G_p\to Q_n$ be reduction over $R_n$, with $Q_n$
defined as its image, and let $Z_n=r_n(Z_p)$. Identify the central
subgroup $Z_n=\{e_{15}(a):a\in R_n\}$ with the additive group of $R_n$.
The character
\[
\eta_n:Z_p\longrightarrow\mathbb T,
\qquad \eta_n(e_{15}(a))=\chi_n(a\bmod(t^{2n+1}-1)),
\]
satisfies
\begin{equation}\label{eq:Thom-characters}
\eta_n(z)\longrightarrow1\quad(z\in Z_p),\qquad
\eta_n(e_{15}(t^n))=\exp(2\pi i/p)\ne1.
\end{equation}
Indeed, for a fixed Laurent polynomial, reduction sends each positive
exponent to a fixed index and each negative exponent $j$ to $2n+1+j$,
once $n$ is sufficiently large. None of these indices is $n$, so the
first limit is eventually equal to $1$.

Let $\rho_n=\operatorname{Ind}_{Z_n}^{Q_n}\chi_n$, of dimension
$d_n=[Q_n:Z_n]$. Since $Z_n$ is central,
\begin{equation}\label{eq:Thom-induced-trace}
\rho_n(z)=\chi_n(z)I_{d_n}\quad(z\in Z_n),\qquad
\tr_{d_n}(\rho_n(q))=
\begin{cases}
\chi_n(q),&q\in Z_n,\\
0,&q\notin Z_n.
\end{cases}
\end{equation}
These identities follow by realizing the induced representation on
the cosets $Q_n/Z_n$. An element outside $Z_n$ permutes these cosets
without fixed points.

Fix a set-theoretic section $s:G_p\to\widetilde G_p$ with $s(e)=e$,
and write
\[
\kappa(g,h)=s(g)s(h)s(gh)^{-1}\in C_p.
\]
The maps
\[
v_n(g)=\rho_n(r_n(s(g))),\qquad
c_n(g,h)=\eta_n(\kappa(g,h))
\]
satisfy
\begin{equation}\label{eq:Thom-projective}
v_n(g)v_n(h)=c_n(g,h)v_n(gh),\qquad
\lim_n c_n(g,h)=1
\quad(g,h\in G_p).
\end{equation}

We now show that none of the multipliers $c_n$ is a coboundary.
Let $a_n,b\in G_p$ be the images of $e_{12}(t^n)$ and $e_{25}(1)$.
They commute, since
\[
[e_{12}(t^n),e_{25}(1)]=e_{15}(t^n)\in C_p.
\]
Their chosen lifts differ from these elementary matrices by elements
of $C_p$, which act as scalars under $\rho_n r_n$. Consequently,
\[
v_n(a_n)v_n(b)=\exp(2\pi i/p)\,v_n(b)v_n(a_n),
\]
and therefore
\[
\frac{c_n(a_n,b)}{c_n(b,a_n)}=\exp(2\pi i/p)\ne1.
\]
A scalar $2$-coboundary is symmetric on every commuting pair. Hence
$c_n$ is not a coboundary for any $n$. This argument applies also
when $p=2$.

We next construct an amenable trace on $C^*(G_p)$.
Set $H_p=\widetilde G_p/Z_p=G_p/N_p$. The reductions induce
homomorphisms $H_p\to Q_n/Z_n$. For each fixed
$x\in\widetilde G_p\setminus Z_p$, at least one entry of its row block,
column block, or middle block minus $I_3$ is a nonzero Laurent
polynomial. This entry remains nonzero modulo $t^{2n+1}-1$ for all
sufficiently large $n$. Thus
\[
r_n(x)\notin Z_n\quad\text{eventually}.
\]
In particular, these finite quotients separate the elements of $H_p$.
The normalized traces of their regular representations, pulled back
to $G_p$, converge pointwise to
\[
\tau(u_g)=\mathbf1_{N_p}(g).
\]
They are traces of genuine finite-dimensional representations of
$C^*(G_p)$, so $\tau$ is amenable by \cite[Theorem~3.1.6]{Brown}.

For $g\in N_p$, its lift $s(g)$ belongs to $Z_p$, and
\eqref{eq:Thom-characters} gives
$\tr_{d_n}(v_n(g))\to1$. For $g\notin N_p$, the preceding separation
argument and \eqref{eq:Thom-induced-trace} give
$\tr_{d_n}(v_n(g))=0$ eventually. Hence
\begin{equation}\label{eq:Thom-trace-limit}
\lim_n\tr_{d_n}(v_n(g))=\tau(u_g)\qquad(g\in G_p).
\end{equation}

Finally, we show that these projective representations obstruct local embeddability.
Fix a free ultrafilter $\omega$ and form the tracial ultraproduct
\[
\mathcal N=\Bigl(\prod_n M_{d_n}\Bigr)/I_\omega,\qquad
I_\omega=\{(x_n):\lim_{n\to\omega}\|x_n\|_{2,\tr_{d_n}}=0\}.
\]
Let $q_\omega$ be its quotient map and $\tau_\omega$ its faithful
normal tracial state. By \eqref{eq:Thom-projective}, there is a unital
$*$-homomorphism
\[
\pi:C^*(G_p)\longrightarrow\mathcal N,\qquad
\pi(u_g)=q_\omega((v_n(g))_n).
\]
Equation~\eqref{eq:Thom-trace-limit} implies
$\tau_\omega\pi=\tau$.

Suppose that $d_f(C^*(G_p))=1$. As $\tau$ is amenable,
Proposition~\ref{prop:trace-LLP} gives local u.c.p. liftings of $\pi$
through $q_\omega$. Property $(T)$ of $G_p$ and
\cite[Theorem~2.3]{ISW} provide a finite set $F\subset G_p$ and
$\delta>0$ such that every projective unitary representation $w$
with a non-coboundary multiplier satisfies
\begin{equation}\label{eq:Thom-projective-gap}
\max_{g\in F}\|w(g)\xi-\xi\|\ge\delta\|\xi\|.
\end{equation}
Lift $\pi$ on the operator system
$E=\spanop\{1,u_g,u_g^*:g\in F\}$. By extending its coordinate maps
using Arveson's theorem, we obtain u.c.p. maps
$\varphi_n:C^*(G_p)\to M_{d_n}$ such that
\[
\lim_{n\to\omega}\|\varphi_n(u_g)-v_n(g)\|_{2,\tr_{d_n}}=0
\qquad(g\in F).
\]
Choose Stinespring dilations
$\varphi_n(a)=V_n^*\sigma_n(a)V_n$, with
$V_n:\ell_2^{d_n}\to H_n$ isometric. On
$\mathcal S_2(\ell_2^{d_n},H_n)$, with its unnormalized
Hilbert--Schmidt norm, set
\[
w_n(g)T=\sigma_n(u_g)T v_n(g)^*,\qquad
\xi_n=d_n^{-1/2}V_n.
\]
Then $\|\xi_n\|=1$, the multiplier of $w_n$ is $\overline{c_n}$,
and
\begin{align*}
\|w_n(g)\xi_n-\xi_n\|^2
&=2-2\Ree\tr_{d_n}\bigl(\varphi_n(u_g)v_n(g)^*\bigr)\\
&\le2\|\varphi_n(u_g)-v_n(g)\|_{2,\tr_{d_n}}.
\end{align*}
The maximum over $g\in F$ tends to zero along $\omega$, contradicting
\eqref{eq:Thom-projective-gap}. Thus $d_f(C^*(G_p))\ne1$, and
Proposition~\ref{prop:df-dichotomy} yields $d_f(C^*(G_p))=\infty$.
\end{proof}

\medskip

\begin{remark}\label{rem:free-group-products}
We recently proved in \cite{Roy26} that other full group $C^*$-algebras
do not locally embed into $C^*(\F_\infty)$, namely
$C^*(\F_k\times\F_m)$ for all integers $k,m\ge2$.
\end{remark}

\medskip

\noindent \textbf{Acknowledgments.} I am grateful to Gilles Pisier for useful and stimulating conversations about various questions on WEP and LLP.\\
\textbf{AI Assistance.} The writing of this paper started in the pre-AI era, but ChatGPT was then used to develop the proofs in the last section.


\begin{thebibliography}{Brown06}

\bibitem[Brown06]{Brown}
N.~P.~Brown,
\emph{Invariant means and finite representation theory of $C^*$-algebras},
Mem. Amer. Math. Soc. \textbf{184} (2006), no.~865.
\href{https://arxiv.org/abs/math/0304009}{arXiv:math/0304009}.

\bibitem[BO08]{BO}
N.~P.~Brown and N.~Ozawa,
\emph{$C^*$-Algebras and Finite-Dimensional Approximations},
Graduate Studies in Mathematics, vol.~88,
American Mathematical Society, Providence, RI, 2008.
\href{https://doi.org/10.1090/gsm/088}{doi:10.1090/gsm/088}.

\bibitem[CE76]{CE}
M.-D.~Choi and E.~G.~Effros,
\emph{The completely positive lifting problem for $C^*$-algebras},
Ann. of Math. (2) \textbf{104} (1976), no.~3, 585--609.
\href{https://doi.org/10.2307/1970968}{doi:10.2307/1970968}.

\bibitem[ES24]{ES}
D.~Enders and T.~Shulman,
\emph{On the (Local) Lifting Property},
\href{https://arxiv.org/abs/2403.12224v2}{arXiv:2403.12224v2} (2024).
The reference to Proposition~5.7 uses this version.

\bibitem[GS15]{GS}
I.~Goldbring and T.~Sinclair,
\emph{On Kirchberg's embedding problem},
J. Funct. Anal. \textbf{269} (2015), no.~1, 155--198.
\href{https://doi.org/10.1016/j.jfa.2015.02.016}{doi:10.1016/j.jfa.2015.02.016}.

\bibitem[Har98]{Harcharras}
A.~Harcharras,
\emph{On some ``stability'' properties of the full $C^*$-algebra associated
to the free group $\F_\infty$},
Proc. Edinb. Math. Soc. \textbf{41} (1998), no.~1, 93--116.
\href{https://doi.org/10.1017/S0013091500019441}{doi:10.1017/S0013091500019441}.

\bibitem[ISW20]{ISW}
A.~Ioana, P.~Spaas, and M.~Wiersma,
\emph{Cohomological obstructions to lifting properties for full
$C^*$-algebras of property $(T)$ groups},
\href{https://arxiv.org/abs/2006.01874v2}{arXiv:2006.01874v2} (2020).
The references to Theorem~2.3 and Sections~3.2--3.3 use this version.

\bibitem[JP95]{JP}
M.~Junge and G.~Pisier,
\emph{Bilinear forms on exact operator spaces and $\mathbb B(H)\otimes \mathbb B(H)$},
Geom. Funct. Anal. \textbf{5} (1995), no.~2, 329--363.
\href{https://doi.org/10.1007/BF01895670}{doi:10.1007/BF01895670};
\href{https://arxiv.org/abs/math/9308208}{arXiv:math/9308208}.

\bibitem[Kir93]{Kirchberg}
E.~Kirchberg,
\emph{On non-semisplit extensions, tensor products and exactness of group
$C^*$-algebras},
Invent. Math. \textbf{112} (1993), 449--489.
\href{https://doi.org/10.1007/BF01232444}{doi:10.1007/BF01232444}.

\bibitem[Pis20]{Pisier}
G.~Pisier,
\emph{Tensor Products of $C^*$-Algebras and Operator Spaces:
The Connes--Kirchberg Problem},
London Mathematical Society Student Texts, vol.~96,
Cambridge University Press, Cambridge, 2020.
\href{https://doi.org/10.1017/9781108782081}{doi:10.1017/9781108782081}.
See also the author's corrected version dated February 22, 2021,
\href{https://www.math.tamu.edu/pisier/TPCOS.pdf}{TPCOS.pdf},
for the numbering of the results cited here.

\bibitem[RS89]{RS}
A.~G.~Robertson and R.~R.~Smith,
\emph{Liftings and extensions of maps on $C^*$-algebras},
J. Operator Theory \textbf{21} (1989), no.~1, 117--131.

\bibitem[Roy26]{Roy26}
J.~Roydor,
\emph{$C^*(\F_2\times\F_2)$ fails the local lifting property},
\href{https://arxiv.org/abs/2611.31147}{arXiv:2611.31147} (2026).

\bibitem[Sin17]{Sinclair}
T.~Sinclair,
\emph{CP-stability and the local lifting property},
New York J. Math. \textbf{23} (2017), 739--747.
\href{https://nyjm.albany.edu/j/2017/23-34.html}{nyjm.albany.edu/j/2017/23-34.html}.

\bibitem[Thi19]{Thiel}
H.~Thiel,
\emph{Inductive limits of projective $C^*$-algebras},
J. Noncommut. Geom. \textbf{13} (2019), no.~4, 1435--1462.
\href{https://arxiv.org/abs/1105.1979}{arXiv:1105.1979}.

\bibitem[Tho10]{Thom}
A.~Thom,
\emph{Examples of hyperlinear groups without factorization property},
Groups Geom. Dyn. \textbf{4} (2010), no.~1, 195--208.
\href{https://doi.org/10.4171/GGD/80}{doi:10.4171/GGD/80}.

\bibitem[Voi90]{Voiculescu}
D.~Voiculescu,
\emph{Property $T$ and approximation of operators},
Bull. London Math. Soc. \textbf{22} (1990), no.~1, 25--30.
\href{https://doi.org/10.1112/blms/22.1.25}{doi:10.1112/blms/22.1.25}.

\end{thebibliography}
\end{document}